\documentclass[11pt]{article}
\usepackage{amsfonts}
\usepackage{amstext}
\usepackage{amsthm}
\usepackage{makeidx}
\usepackage{amsmath}
\usepackage{eufrak}
\usepackage{wrapfig}
\usepackage{mathrsfs}
\usepackage{amssymb}
\usepackage{latexsym}
\usepackage[latin1]{inputenc}
\usepackage[english]{babel}
\usepackage{graphicx}
\input{xypic.tex} \xyoption{all}

\newtheorem{teo}{Theorem}[section]

\newtheorem{ddef}[teo]{Definition}
\newtheorem{ex}[teo]{Example}
\newtheorem{cor}[teo]{Corollary}

\newtheorem{obs}{Remark}
\newtheorem{lem}[teo]{Lemma}

\newcommand{\F}{{\mathcal{F}}}
\newcommand{\wF}{\widetilde{{\mathcal{F}}}}
\newcommand{\G}{{\mathcal{G}}}
\newcommand{\wG}{\widetilde{{\mathcal{G}}}}
\newcommand{\D}{{\mathcal{D}}}
\newcommand{\wD}{\widetilde{{\mathcal{D}}}}
\newcommand{\pn}{{\mathbb{P}^n}}
\newcommand{\C}{{\mathbb{C}}}

\begin{document}
\title{\textsc{inequalities for characteristic numbers of flags of distributions and foliations}}
\author{Maurício Corrêa Jr. and Márcio G. Soares}
\date{}
\maketitle
\begin{abstract}
\noindent We prove inequalities relating the degrees of holomorphic distributions and of holomorphic foliations forming a flag on $\mathbb{P}^n$. Such inequalities are inspired by the so called Poincar\'e problem for foliations. \end{abstract}

 \footnotetext[3]{{\sl This work was partially supported by CNPq, CAPES and FAPEMIG.} }

\section{Introduction}

In this paper we consider flags of distributions and of foliations on complex projective spaces and deduce inequalities relating their degrees.

Before stating the results we recall that a holomorphic \textsl{distribution}, or a \textsl{Pfaff equation}, on a complex manifold $M$, is defined by a holomorphic line bundle $\mathcal{L}$ on $M$ and a nontrivial global section $\omega \in H^0(M, \Omega^p_M \otimes \mathcal{L})$, where $\Omega^p_M$ is the sheaf of holomorphic $p$-forms on $M$. The number $p$, $1 \leq p \leq n-1$, is the codimension of the distribution, where $n=\dim M$. A holomorphic \textsl{foliation} is obtained by imposing the \textsl{integrability} condition to a distribution and, this being the case, the line bundle $\mathcal{L}$ corresponds to the determinant bundle of the rank $p$ normal sheaf of the foliation.

To a distribution on $\mathbb{P}^n$, and hence to a foliation, we can associate a nonnegative integer, its \textsl{degree}, which is the degree of the variety formed by the points $x \in \mathbb{L}^p$, a fixed generic linear subspace of dimension $p$, at which the $(n-p)$-plane of the distribution, passing through the point $x$, is not in general position with respect to this subspace. Also, by a \textsl{flag} of distributions, $\mathscr{D}\,:= (\D_{j_1}, \D_{j_2}, \dots, \D_{j_m})$, we mean a collection of distributions of dimensions $1 \leq j_1 < j_2 < \cdots < j_m < n$ such that, at each point $x$ where the distributions are regular, $D_{{j_r},x} \subset D_{{j_s},x}$ whenever $r < s$. All these notions are explained in Section \ref{defs}.

The results are

\begin{teo}\label{dist}  Let $\mathscr{D}\,:=(\F, \mathcal{G})$ be  a flag of reduced holomorphic distributions on  $\mathbb{P}^n$, $n \geq3$, with  $\dim(\F)=\mathrm{codim}(\mathcal{G})=1$ and $\deg(\G) \geq 2$. If
\begin{itemize}
\item [{(i)}] $\deg(\F)\neq \displaystyle\left(\frac{n}{2}\right) \deg(\G)$ if $n$ is even;
\item [(ii)] $\deg(\F)\neq \displaystyle\left(\frac{n-1}{2}\right) \deg(\G)-1$ if $n$ is odd;
\item[(iii)] $\mathrm{Sing}(\G)$ is isolated;
\end{itemize}
then $\deg(\mathcal{G})\leq\deg(\F)-1$.
\end{teo}

\begin{obs}
Since $\deg(\G) \geq 2$ we always have $(\frac{n}{2}) \deg(\G)\neq \deg(\G)+1$ and $(\frac{n-1}{2}) \deg(\G)\neq \deg(\G)+1$.
\end{obs}

\begin{teo}\label{dist1}  Suppose $\mathscr{F}\,:=(\F,\G)$ is a flag of reduced holomorphic distributions on $\mathbb{P}^n$, with $\dim(\F)=k $, $\mathrm{codim}(\G)=1$ and $\mathrm{codim} (\mathrm{Sing} (\G))\geq n-k+1$. If  the tangent sheaf $\wF$ of $\F$ is split,  then $\deg(\mathcal{G})\leq \deg(\F)$.
\end{teo}

\begin{teo}\label{fol1}  Let $\mathscr{F}\,:=(\F,\G)$ be a flag of reduced holomorphic foliations on $\mathbb{P}^n$, $n\geq3$.
If $\dim (\F) = \dim (\G)-1$ and   $\mathrm{Sing} (\G)$ has a  Baum-Kupka component ${K}$, then
$$
\deg(\mathcal{G})\leq \deg(\F).
$$
\end{teo}

\begin{cor}\label{flag}
Let $\mathscr{F}\,:= (\F_1, \F_2, \dots, \F_k)$ be a flag of reduced foliations on $\mathbb{P}^n$ with $\dim (\F_j)= \dim (\F_{j+1}) -1$ for
$j= 1, 2, \dots, k-1$. If $\mathrm{Sing} (\F_{j+1})$ has a Baum-Kupka component $K_{j+1}$, for $j= 1, 2, \dots, k-1$, then
$$
\deg(\F_1) \leq \deg(\F_2) \leq \dots \leq \deg(\F_k).
$$
\end{cor}

\section{Preliminaries}\label{defs}

We start by recalling some definitions.

\begin{ddef}\label{fol}
Let $M$ be a connected complex manifold of dimension $n$ and $ \mathcal{O}(TM)$ be its tangent sheaf. A singular holomorphic distribution $\D$ on $M$, of dimension $r$, is a coherent subsheaf $\wD$ of $\mathcal{O}(TM)$ of rank $r$. In case $\wD$ is involutive (or integrable) we have a singular holomorphic foliation on $M$, of dimension  $r$. Integrable means that, for each $p \in M$, the stalk $\wD_p$ is closed under the Lie bracket operation, $[\wD_p, \wD_p \,]\subset \wD_p$.
\end{ddef}

In the above, the rank of $\wD$ is the rank of its locally free part. Since $\mathcal{O}(TM)$ is locally free, the coherence of $\wD$ simply means that it is locally finitely generated. We call
$\wD$ the {\it tangent sheaf} of the distribution and the quotient,
$\mathcal{N}_{\D}=\mathcal{O}(TM)/\wD$, its {\it normal sheaf}.

The \emph{singular set} of $\D$ is defined by
$$
S(\D) = \{ p \in M \,:\, (\mathcal{N}_{\D})_p \; \mathrm{is\;
not\; a\; free \;} \mathcal{O}_p-\mathrm{module}\}.
$$
In case we have a foliation we will use the notation $\F$, for the foliation, and $\wF$ for its tangent sheaf. On $M \setminus
S(\F)$ there is a unique (up to isomorphism) holomorphic vector subbundle $E$ of
the restriction ${TM}_{\vert\,{M \setminus S(\F)}}$, whose sheaf of
germs of holomorphic sections, $\widetilde{E}$, satisfies
$\widetilde{E}= \wF_{\vert\,{M \setminus S(\F)}}$. Clearly $r=$ rank of $E$.

We will assume that $\wD$ is {\sl full} (or saturated) which
means: let $U$ be an open subset of $M$ and $\sigma$ a holomorphic section of $\mathcal{O}(TM)_{|U}$ such that $\sigma_p \in \wD_p$ for all $p
\in U \cap (M \setminus S(\D))$. Then we have that for all $p \in U$,
$\sigma_p \in\wD_p$. In this case the distribution (or foliation, if this is the case) $\D$ is said to be {\sl reduced}.

An equivalent formulation of \textsl{full} is as follows: let $\Omega^1 = \mathcal{O}(T^\ast M)$ be the cotangent sheaf of $M$. Set $\wD^o = \{ \omega \in
\Omega^1\,:\, i_{\gamma}\omega =0\,\,\forall \;\gamma \in \wD\}$ and
$\wD^{oo} = \{ \gamma \in \mathcal{O}(TM)\,:\, i_{\gamma}\omega
=0\,\,\forall \;\omega  \in \wD^o\}$, where $i$ is the contraction. $\wD$ is full if
$\wD = \wD^{oo}$. Note that integrability of
$\wD$ implies integrability of $\wD^{oo}$.

Singular distributions and foliations can dually be defined in terms of the cotangent
sheaf. Thus a {\sl singular distribution of corank $q$}, $\G$, is a coherent
subsheaf $\wG$ of rank $q$ of $\Omega^1$. $\wG$ is called the {\sl conormal} sheaf of the distribution $\D$. Its annihilator
\[
\D=\G^o=\{\, \gamma \in \mathcal{O}(TM)\,:\, i_\gamma \omega=0\ \ \text{for all}
\ \omega \in\wG\,\}
\]
 is a singular distribution of rank $r=n-q$. The singular set of $\G$, $\mathrm{Sing}(\G)$, is the set $\mathrm{Sing}(\Omega^1 / \wG)$. See T. Suwa \cite{Suw} for the relation between these two definitions.

We remark that, if a foliation $\F$ is reduced then $\mathrm{codim}\, S(\F) \geq 2$ and reciprocally, provided $\wF$ is locally free (see \cite{Suw}). This is a useful concept since it avoids the appearance of ``fake" (or ``removable") singularities.

\begin{ddef}\label{flag}

Let $\D_{j_1}, \D_{j_2}, \dots, \D_{j_m}$ be holomorphic distributions (foliations) on a connected complex manifold $M^n$. They form a $\textrm{flag}$ provided\\

    (i) $1\leq j_1 < j_2 < \cdots < j_m < n= \dim M$ and $\dim \D_{j_i} = j_i$.\\

    (ii) $\wD_{j_i}$ is a subsheaf of $\wD_{j_{i+1}}$. Here, $\wD_{j_r}$ is the tangent sheaf of $\D_{j_r}$.
\end{ddef}

\begin{obs}\label{obsflag} For foliations, outside $\mathrm{Sing}(\F_{j_i}) \cup \mathrm{Sing}(\F_{j_{r}})$, $j_i <j_{r}$, we have $T_p \F_{j_i} \subset T_p \F_{j_r}$, so that the leaves of $T_p \F_{j_r}$ are foliated by the leaves of $T_p \F_{j_i}$. By a result of J. Yoshizaki \cite{Y} (see also R. Mol \cite{Mol}) the singular set $\mathrm{Sing}(\F_{j_{r}})$ is invariant by $\F_{j_i}$ whenever $j_i <j_{r}$.
\end{obs}

As for the structure of the singular set of a foliation of dimension $r$ we have the following result
of P.Baum \cite{Baum}, in the version due to J.B.Carrell \cite{Carrell} in
the review of \cite{Baum} (this result also appears in \cite{Cerveau}):

\begin{teo}\label{debaum}
Let $p$ be a smooth point of $\mathrm{Sing}(\F)$ with  $\dim T_p
\mathrm{Sing}(\F) = \dim \wF(p) = r-1$, where $\wF(p)=\{v(p)\; | \; v \in \wF_p\}$. Then there exists a
neighborhood $U_p \subset M$ of $p$ and a holomorphic submersion
$f: U_p \to {\mathbb{C}}^{n-r+1}$, $f(p)=0$, such that $f^{-1}(0) = U_p
\cap \mathrm{Sing}(\F)$ and such that $\wF_{\vert {U_p}}= (f^{\ast}
\xi^o)^o$, where $\xi$ is the sheaf on $f(U_p)$ generated by a
holomorphic vector field $X$ on $f(U_p)$ with its only zero at
$0$.
\end{teo}

It follows that the foliation $\F$ is, in $U_p$, the pull-back
via $f$ of the foliation $\widehat{\F}$ induced by $X$ in
$f(U_p)$ and, hence, we have a local product structure. We call
such singularities of {\sl Baum-Kupka} type in view of a prior result of I. Kupka \cite{Kupka}
for codimension one holomorphic foliations which states that, if $\F$ is
given by the integrable one-form $\omega$ and $p$ is a point such
that $\omega(p) =0$ and $d\,\omega(p) \neq 0$ then, in a neighborhood
of $p$, $\F$ is the pull-back via a submersion of a
one-dimensional foliation defined around $0 \in {\mathbb{C}}^2$ and with
an isolated singularity at $0$.

\subsection{The case of $\pn$}

\begin{ddef}\label{degree}
Let $\D$ be a codimension $n-k$ distribution on $\mathbb P^n$ given
by $\omega \in H^0(\mathbb P^n, \Omega^{n-k}_{\mathbb P^n} \otimes
\mathcal L)$. If $\mathrm{i}: \mathbb P^{n-k} \to \mathbb P^n$ is a general
linear immersion then $\mathrm{i}^* \omega \in H^0(\mathbb P^{n-k},
\Omega^{n-k}_{\mathbb P^{n-k}} \otimes \mathcal L)$ is a section of a line
bundle, and its zero divisor reflects the tangencies between
$\D$ and $\mathrm{i}(\mathbb P^{n-k})$. The degree of $\D$ is the degree of such tangency divisor. It is noted $\deg(\D)$.
\end{ddef}

Set $d:=\deg(\D)$. Since $\Omega^{n-k}_{\mathbb P^{n-k}}\otimes
\mathcal L= \mathcal O_{\mathbb P^{n-k}}( \deg(\mathcal L) - n+k - 1)$, one concludes that
$\mathcal L= \mathcal O_{\mathbb P^n}(d+ n-k + 1)$. Besides, the Euler sequence implies that a section $\omega$ of
$\Omega^{n-k}_{\mathbb P^n} ( d + n-k + 1  )$ can be
thought as a polynomial $(n-k)$-form on $\C^{n+1}$ with  homogeneous
coefficients of degree $d + 1$, which we will still
denote by $\omega$, satisfying
\begin{equation}
\label{equirw}
i_\vartheta  \omega = 0
\end{equation}
where $\vartheta=x_0 \frac{\partial}{\partial x_0} + \cdots + x_n
\frac{\partial}{\partial x_n}$ is the radial vector field and $i_\vartheta$ means contraction by $\vartheta$. Thus the
study of distributions of degree $d$ on $\mathbb P^n$ reduces to the
study of locally decomposable homogeneous $(n-k)$-forms  on $\mathbb C^{n+1}$, of
degree $d+1$, satisfying relation (\ref{equirw}).

Let $\wD$ be the tangent sheaf of $\D$. If the singular set of $\D$ has
codimension at least two we obtain the adjunction formula
$$
K_{\pn}=\det(\wD^*)\otimes \det(\mathcal{N}_{\D}^*).
$$
Since $\det(\mathcal{N}_{\D}^*)= \mathcal O_{\mathbb P^n}(-d- n+k - 1)$ and $K_{\pn}=\mathcal O_{\mathbb P^n}(-n-1)$, then $\det(\wD)=\mathcal O_{\mathbb P^n}(k-d)$.

We close with a definition. This is motivated by the fact that the singular set of a codimension one foliation on $\pn$ has at least
a codimension two irreducible component.

\begin{ddef}\label{BKcomp}
Let $\F$ be a foliation on $\pn$, $n \geq 3$, of codimension $n-k$. An analytic subset $K \subset \mathrm{Sing}(\F)$, of codimension
$n-k+1$, is a Baum-Kupka component if $K$ is an irreducible component of $\mathrm{Sing}(\F)$ whose points are all singularities of
Baum-Kupka type and, moreover, if $\omega \in H^0(\mathbb P^n, \Omega^{n-k}_{\mathbb P^n}(d+n-k+1))$ induces $\F$, then $d \omega_{|K}$ is nowhere vanishing.
\end{ddef}

\section{Proofs}

\subsection{Proof of Theorem \ref{dist}}

We will make use of Bott's formulae:

Let $n,p,q$ and $k$ be integers,
with $n$ positive and $p$ and $q$ nonnegative.  Then
\begin{equation*}
h^q(\mathbb P^n,\Omega_{\mathbb P^n}^p(k)) =
\begin{cases}
\binom{k+n-p}{k}\binom{k-1}{p} & \text{for } q=0, 0\le p\le n \text{ and } k>p,\\
1 & \text{for } k=0 \text{ and } 0\le p=q\le n,\\
\binom{-k+p}{-k}\binom{-k-1}{n-p} & \text{for } q=n, 0\le p\le n \text{ and } k<p-n,\\
0 & \text{otherwise.}
\end{cases}
\end{equation*}
Let $r,s$ and $t$ be integers, with $r$ and $s$ nonnegative. Observe
that the natural pairing $\Omega^p_{\mathbb P^n}\otimes\Omega^{n-p}_{\mathbb P^n}
\to \Omega^n_{\mathbb P^n}$ is perfect and so there is an induced isomorphism
$\wedge^{r}T_{\mathbb P^n}(t)\simeq \Omega_{\mathbb P^n}^{n-r}(t+n+1)$. Hence the
formulae above become
\begin{equation*}
h^s(\mathbb P^n,\wedge^{r}T_{\mathbb P^n}(t)) =
\begin{cases}
\binom{t+n+1+r}{t+n+1}\binom{t+n}{n-r} & \text{for } s=0, 0\le r \le n \text{ and } t+r\ge 0,\\
1 & \text{for } t=-n-1 \text{ and } 0\le n-r=s\le n,\\
\binom{-t-1+r}{-t-n-1}\binom{-t-n-2}{r} & \text{for } s=n, 0\le r\le n \text{ and } t+n+r+2\le 0,\\
0 & \text{otherwise.}
\end{cases}
\end{equation*}

Now for the proof.

Let $\omega\in H^0(\mathbb P^n,\Omega^1_{\mathbb P^n}(m+2))$ be a global section
inducing the distribution $\G$ of degree $m$. Suppose there exists a one-dimensional reduced distribution $\F$ (necessarily a foliation) of degree $d$, tangent to $\G$ and induced by $X\in H^0(\mathbb P^n,T{\mathbb P^n}(d-1))$.
$\omega$ induces a morphism of sheaves
\begin{equation}\label{B1}
T{\mathbb P^n}(d-1)\stackrel{\imath_{\omega}} \longrightarrow  \mathcal{O}(d+m+1).
\end{equation}
By hypothesis, since the set of zeros $Z$ of $\omega\in
H^0(\mathbb P^n,\Omega^1_{\mathbb P^n}(m+2))$ is  isolated, we have that Koszul
complex associated to $Z$,
\begin{equation}\label{B2}
0 \to \bigwedge^n (\Omega^1_{\mathbb P^n}(m+2))^* \to  \cdots \to
\bigwedge^{2}(\Omega^1_{\mathbb P^n}(m+2))^* \to  (\Omega^1_{\mathbb P^n}(m+2))^*
\to I_Z \to 0,
\end{equation}
is exact.

Rewriting
\begin{equation}\label{B3}
0 \to \bigwedge^n T\mathbb P^n(-n(m+2)) \to  \cdots \to
\bigwedge^{2}T\mathbb P^n(-2m-4) \to  T\mathbb P^n(-m-2) \to I_Z \to
0
\end{equation}
and tensorizing by $\mathcal{O}(d+m+1)$ we obtain the exact sequence
\begin{equation}\label{B4}
\cdots \to \bigwedge^{2}T\mathbb P^n(d-m-3) \to  T\mathbb P^n(d-1)
\stackrel{\imath_{\omega}} \longrightarrow
I_Z(d+m+1) \to 0.
\end{equation}

The sequence in (\ref{B4}) breaks into short exact sequences:
\begin{equation}\label{B5}
\begin{array}{c}
0 \to\bigwedge^nT_{\mathbb P^n}(\underbrace{-n(m+2)+d+m+1}\limits_{\displaystyle= t_n})\to \bigwedge^{n-1}T_{\mathbb P^n}(\underbrace{-(n-1)(m+2)+d+m+1}\limits_{\displaystyle= t_{n-1}})\\
\to K_{n-2}\to 0 \hfill\\
\vdots \\
0 \to K_{r} \to  \bigwedge^{r}T_{\mathbb P^n}(\underbrace{-r(m+2)+d+m+1}\limits_{\displaystyle = t_r}) \to K_{r-1} \to 0\\
\vdots\\
0 \to K_{2} \to  \bigwedge^{2}T_{\mathbb P^n}(\underbrace{d-m-3}\limits_{\displaystyle = t_2}) \to K_{1} \to 0\\
 \\
0 \to K_{1} \to   T\mathbb P^n(\underbrace{d-1}\limits_{\displaystyle = t_1})
\stackrel{\imath_{\omega}} \longrightarrow
 \mathcal{O}(d+m+1) \to 0.
\end{array}
\end{equation}

We have that $X\in H^0(\mathbb P^n,T{\mathbb P^n}(d-1))$ inducing $\F$ satisfies $X\in H^0(\mathbb P^n, K_{1})$.
Passing to cohomology we get the  exact sequence
\begin{equation}\label{B6}
H^0(\mathbb P^n,\bigwedge^{2}T_{\mathbb P^n}(d-m-3) ) \stackrel{\zeta} \longrightarrow
H^0(\mathbb P^n,K_1)\longrightarrow H^1(\mathbb P^n, K_2).
\end{equation}
If we can prove that $\zeta$ is surjective, then it will exist a nonzero global section $\eta$ of
$\bigwedge^{2}T\mathbb P^n(d-m-3)$ such that $\zeta(\eta)=X$.

In particular $ H^0(\mathbb P^n,\bigwedge^{2}T\mathbb P^n(d-m-3))\neq 0$. But, by
Bott's formulae this is possible only when $(d-m-3)+2\geq 0$, and theorem is proved. $\hfill \square$

To prove the surjectivity of $\zeta$ it's enough, by (\ref{B6}), to show that $H^1(\mathbb P^n,K_2)=0$.
This is done in the
\begin{lem}\label{tech} $H^1(\mathbb P^n,K_2)=0$ whenever $d \neq (\frac{n}{2})\,m $ if $n$ is even and $d \neq (\frac{n-1}{2})\,m-1$ if $n$ is odd.
\end{lem}

\noindent\textbf{Proof}  The short sequences in (\ref{B5}) give rise to the cohomology sequences
\begin{equation}\label{B8}
\begin{array}{c}
\dots \rightarrow H^{n-3}(\mathbb P^n,\bigwedge^{n-1}T_{\mathbb P^n}(t_{n-1}) )   \rightarrow
H^{n-3}(\mathbb P^n,K_{n-2})\to H^{n-2}(\mathbb P^n,\bigwedge^{n}T_{\mathbb P^n}(t_n) )\to\\
\\
\to H^{n-2}(\mathbb P^n,\bigwedge^{n-1}T_{\mathbb P^n}(t_{n-1}) )\to\dots
\end{array}
\end{equation}
and, for $3\leq r \leq n-2$,
\begin{equation}\label{B8j}
\begin{array}{c}
\dots \to H^{r-2}(\mathbb P^n,\bigwedge^{r}T_{\mathbb P^n}(t_r) )   \rightarrow
H^{r-2}(\mathbb P^n,K_{r-1})\to H^{r-1}(\mathbb P^n, K_r)\to\\
\\
H^{r-1}(\mathbb P^n,\bigwedge^{r}T_{\mathbb P^n}(t_r) )\to \dots
\end{array}
\end{equation}

Let us first analyze (\ref{B8}).\\

\textsl{Case (i)} $n=3$. This gives $K_2= \bigwedge^{3}T_{\mathbb P^3}(t_3)$ and Bott's formulae tells us that
$H^{1}(\mathbb P^3,\bigwedge^{3}T_{\mathbb P^3}(t_3) )=0$. In this case (\ref{B6}) and (\ref{B8}) coincide and Theorem \ref{dist} holds.

\textsl{Case (ii)} $n >3$. Here we calculate $H^{n-3}(\mathbb P^n,\bigwedge^{n-1}T_{\mathbb P^n}(t_{n-1}) ) $ and \break
$H^{n-2}(\mathbb P^n,\bigwedge^{n-1}T_{\mathbb P^n}(t_{n-1}) )$.

For $s= n-3$ and $r=n-1$ the only possibilities to have \break
$H^{n-3}(\mathbb P^n,\bigwedge^{n-1}T_{\mathbb P^n}(t_{n-1}) ) \neq 0$ are

(I) $s=n-3=0$, $0 \leq r \leq n$ and $t_{n-1} +r \geq 0$, which is ruled out since $n > 3$.

(II) $t_{n-1} = -n-1$ and $0 \leq n-r =s \leq n$. Now, $s=n-3=n-r=n-(n-1)=1$ gives $n=4$. On the other hand,
$-n-1=t_{n-1}= -(n-1)(m+2)+d+m+1$ gives $d=2m=(\frac{4}{2})m$, which is contrary to the hypotheses.

For $s= n-2$ and $r=n-1$ the only possibility to have \break
$H^{n-2}(\mathbb P^n,\bigwedge^{n-1}T_{\mathbb P^n}(t_{n-1}) ) \neq 0$ is $t_{n-1} = -n-1$ and $0 \leq n-r =s \leq n$.
But this gives $s=n-2=n-r=n-(n-1)=1$ and then $n=3$, contradicting $n>3$.

Hence, for $n>3$,
\begin{equation}\label{vanish1}
\begin{array}{c}
H^{n-3}(\mathbb P^n,\bigwedge^{n-1}T_{\mathbb P^n}(t_{n-1}) ) =0\\
\\
H^{n-2}(\mathbb P^n,\bigwedge^{n-1}T_{\mathbb P^n}(t_{n-1}) )=0
\end{array}
\end{equation}

and from (\ref{B8}) we get
\begin{equation}\label{iso1}
H^{n-3}(\mathbb P^n,K_{n-2})\simeq H^{n-2}(\mathbb P^n,\bigwedge^{n}T_{\mathbb P^n}(t_n) ).
\end{equation}

Let us now analyze the sequences in (\ref{B8j}).\\

For $r=3$, by Bott's formulae,
\begin{equation}\label{B9}
H^1(\mathbb P^n,\bigwedge^{3}T_{\mathbb P^n}(t_3) )=H^2(\mathbb P^n,\bigwedge^{3}T_{\mathbb P^n}(t_3) )=0.
\end{equation}
In fact, $H^1(\mathbb P^n,\bigwedge^{3}T_{\mathbb P^n}(t_3) )\neq 0$ if $n-3=1$ and $t_3=-n-1$. This implies that $n=4$ and $d-2m-5=-5$, i.e, $d=2m=(\frac{4}{2})m$, which is ruled out  by hypothesis.

Also,  $H^2(\mathbb P^n,\bigwedge^{3}T_{\mathbb P^n}(t_3) )\neq 0$ only if $n-3=2$ and $t_3=-n-1$. This implies that $n=5$ and $d-2m-5=-6$, i.e, $d=2m-1=(\frac{5-1}{2})m-1$, which is forbidden.
By (\ref{B8j}) we are left with
\begin{equation}\label{B10}
H^1(\mathbb P^n,K_2)\simeq H^2(\mathbb P^n,K_3).
\end{equation}
Repeating this argument we get  the vanishing of
\begin{equation}\label{B11}
\begin{array}{c}
H^{r-2}(\mathbb P^n,\bigwedge^{r}T_{\mathbb P^n}(\underbrace{-r(m+2)+d+m+1}\limits_{\displaystyle = t_r})) =\\
 \\
= H^{r-1}(\mathbb P^n,\bigwedge^{r}T_{\mathbb P^n}(-r(m+2)+d+m+1))=0
\end{array}
\end{equation}
for  $3\leq r\leq n-2$.

Indeed, if $1\leq r-2\leq n-3$, by Bott's formulae,  $H^{r-2}(\mathbb P^n,\bigwedge^{r}T_{\mathbb P^n}(t_r))\neq 0$ only if $n-r=r-2$ and $-r(m+2)+d+m+1=-n-1$. This amounts to $d=(\frac{n}{2})m$, which is not allowed.
The case $H^{r-1}(\mathbb P^n,\bigwedge^{r}T_{\mathbb P^n}(-r(m+2)+d+m+1))\neq 0$ is dealt with analogously and we arrive at $d = (\frac{n-1}{2})\,m-1$, which is forbidden.

From (\ref{B8j}), (\ref{B10}) and (\ref{B11}) we have
\begin{equation}\label{B12}
H^1(\mathbb P^n,K_2)\simeq H^2(\mathbb P^n,K_3)\simeq \cdots \simeq  H^{n-3}(\mathbb P^n,K_{n-2}).
\end{equation}

Invoking (\ref{iso1}) and (\ref{B12}) we conclude that
\begin{equation}\label{B15}
H^1(\mathbb P^n,K_2)\simeq  H^{n-2}(\mathbb P^n,\bigwedge^{n}T_{\mathbb P^n}(t_n) )
\end{equation}
By Bott's formulae $ H^{n-2}(\mathbb P^n,\bigwedge^{n}T_{\mathbb P^n}(t_n) )=0$ and then
$H^1(\mathbb P^n,K_2)=0$. This finishes the proof of the Lemma and of Theorem \ref{dist}.
$\hfill \square$

\subsection{Proof of Theorem \ref{dist1}}

Let $dV =dz_0\wedge \cdots \wedge dz_n$ and $X_1,\dots,X_k$ be homogeneous vector fields
such that $T\F=\bigoplus T\F_{X_i} $. Thus, $\deg(\F)=\sum_{i=1}^k\deg(X_i)$.
Consider the $(n-k)$-form
$i_{X_1}\cdots i_{X_k} i_{\vartheta} dV$, where $\vartheta= \sum\limits_{0}^n z_i\, \frac{\partial}{\partial z_i}$ is the radial vector
field and $i_Y \eta$ is the contraction of $\eta$ by $Y$.

Suppose $\G$ is given
by $\omega \in H^0(\mathbb P^n, \Omega^{1}_{\mathbb P^n} \otimes
\mathcal O_{\mathbb P^n}(\deg(\G)+ 2))$. Then
\begin{equation}\label{contracao}
(i_{X_1}\cdots i_{X_k}i_{\vartheta}dV)\wedge \omega=0.
\end{equation}
In fact, since $i_{X_1}\omega =\dots = i_{X_k}\omega= i_{\vartheta}\omega=0$, we
have
$$
0=i_{\vartheta}(dV\wedge \omega)= (i_{\vartheta}dV)\wedge
\omega+(-1)^{n+1}dV\wedge( i_{\vartheta} \omega)=(i_{\vartheta}dV)\wedge \omega
$$
and
$$
0=i_{X_k}[(i_{\vartheta}dV)\wedge \omega]=(i_{X_k}i_{\vartheta}dV)\wedge
\omega+(-1)^n(i_{\vartheta}dV)\wedge (i_{X_k}
\omega)=(i_{X_k}i_{\vartheta}dV)\wedge \omega.
$$
Proceeding inductively we obtain (\ref{contracao}).
Now, $\mathrm{codim}( \mathrm{Sing} (\G)) \geq n-k+1$ and $(i_{X_1}\cdots i_{X_k}i_{\vartheta}dV)\wedge \omega=0.$
This allow us to invoke Saito's generalization of the de Rham
division Lemma \cite{Saito} and conclude that there exists a homogeneous polynomial
$(n-k-1)$-form $\eta$ on $\mathbb{C}^{n+1}$ such that
$$i_{X_1}\cdots i_{X_k}i_{\vartheta}dV= \omega \wedge \eta.$$
Computing degrees,
$$
\begin{array}{c}
\deg(\F)+1=\sum\limits_{i=1}^k\deg(X_i)+1=\deg(i_{X_1}\cdots
i_{X_k}i_{\vartheta}dZ)  \\
\\
= \deg(\omega \wedge
\eta)=\deg(\omega)+\deg(\eta)=\deg(\G)+1+\deg(\eta)
\end{array}
$$
and thus $\deg(\G)\leq \deg(\F)$.

$\hfill \square$

\begin{ex}\label{disti} A codimension one distribution.
\end{ex}
A generic codimension one distribution on $\pn$, of degree $k$, has as singular locus a zero dimensional smooth algebraic variety of degree $\displaystyle \frac{(k+1)^{n+1} - (-1)^{n+1}}{k+2}$ (see \cite{J}, Th. 2.3, pg. 87).

Here we show that the bound given in Theorem \ref{dist} $(i)$ is sharp. This example can easily be generalized to any dimension, but we will give it in $\mathbb{P}^3$. Consider the antisymmetric matrix
$$
M = \left(
  \begin{array}{cccc}
    0 & 0 & 0 & z_3^k \\
    0 & 0 & z_2^k & z_0^k \\
    0 & -z_2^k & 0 & 0 \\
    -z_3^k & -z_0^k & 0 & 0 \\
  \end{array}
\right)
$$
and let $\omega$ be the 1-form $\omega= \sum\limits_0^3 A_i\;dz_i$ where
$$
\left(
  \begin{array}{c}
    A_0 \\
    A_1 \\
    A_2 \\
    A_4 \\
  \end{array}
\right) =
M
\left(
  \begin{array}{c}
    z_0 \\
    z_1 \\
    z_2 \\
    z_3 \\
  \end{array}
\right)
$$
We have $\sum\limits_0^3 z_i\,A_i \equiv 0$ because $M$ is antisymmetric, so $\omega$ defines a distribution $\D_\omega$ on $\mathbb{P}^3$. As
$$
\omega = z_3^{k+1}\,dz_0 + (z_2^{k+1} + z_0^k z_3)\,dz_1 -z_1 z_2^{k}\, dz_2 +
(-z_0 z_3^k - z_0^k z_1)\, dz_3$$
we have $\deg(\D_\omega)=k$ and $\mathrm{Sing}(\D_\omega) = \{(1:0:0:0), (0:1:0:0)\}$ not counting multiplicities. On the other hand, the foliation $\F$ on $\mathbb{P}^3$, of degree $k+1$, induced by the vector field
$$
X= z_1 z_2^{k} \frac{\partial}{\partial z_0} +(z_0 z_3^k + z_0^k z_1) \frac{\partial}{\partial z_1} + z_3^{k+1} \frac{\partial}{\partial z_2} +
(z_2^{k+1} + z_0^k z_3) \frac{\partial}{\partial z_3}
$$
is tangent to $\D_\omega$ and $\deg(\D_\omega)= \deg(\F)-1$.

\subsection{Proof of Theorem \ref{fol1}}

Since ${K}$ is a Baum-Kupka component of $\mathrm{Sing}(\G)$ we have
$k:=\dim {K}=\dim(\G)-1=\dim(\F)$. We claim that

\begin{equation}\label{canonico}
\Omega_{{K}}^{k}=\mathcal{O}_{{K}}(\deg(\G)-\dim\, (\F)-1).
\end{equation}
To see this, if $\omega\in H^0\left(\pn, \Omega_{\pn}^{n-\dim(\G)}\otimes \mathcal{O}_{\pn}(\deg(\G)+n-\dim(\G)+1)\right)$ is a $(n-\dim(\G))$-form inducing $\G$, then $d \omega|_{{K}}$ defines a nowhere vanishing holomorphic section of $\bigwedge^{n-k} \nu_{{K}}^*\otimes \mathcal{O}_{\pn}(\deg(\G)+n-\dim(\G)+1)_{|{{K}}}$, where $\nu_K$ is the normal sheaf of $K$. In particular,
 $$
 \bigwedge^{n-k} \nu_{{K}}^*\otimes \mathcal{O}_{\pn}\left(\deg(\G)+n-\dim(\G)+1\right)_{|{{K}}}= \!\!\bigwedge^{n-k} \nu_{{K}}^*\otimes \mathcal{O}_{{K}}\left(\deg(\G)+n-\dim(\G)+1\right)
 $$
 is trivial and thus $\bigwedge^{n-k} \nu_{{K}}\simeq  \mathcal{O}_{{K}}(\deg(\G)+n-\dim(\G)+1).$

Now, using the adjunction formula
\begin{center}
$
\Omega_{{K}}^{k}=\Omega_{\pn}^{n}|_{{K}}\otimes \bigwedge^{n-k} \nu_{{K}},
$
\end{center}
$\Omega_{\pn}^{n}=\mathcal{O}_{\pn}(-n-1)$ and $\dim(\G)=\dim(\F)+1$ we conclude (\ref{canonico}).

The foliation $\F$ induces a map  $ \det(\wF) \longrightarrow \bigwedge^k T\pn$, that furnishes a holomorphic global section of
$$
\bigwedge^k T\pn \otimes \det(\wF)^*=\bigwedge^k T\pn \otimes \mathcal{O}_{\pn}(\deg(\F)-k),
$$
because  $\det(\wF)^*=\mathcal{O}_{\pn}(\deg(\F)-k)$.

Since ${K}$ is invariant by $\F$ and ${K} \not\subset \mathrm{Sing} (\F)$, we have that $\F_{|{K}}$ induces a  nonzero  holomorphic global section $\zeta$ of
$$
\bigwedge^k T{K}\otimes\det(\wF)^*_{|{{K}}}=
(\Omega_{{K}}^{k})^{*}\otimes\mathcal{O}_{{K}}(\deg(\F)-k).
$$
It follows from \cite[Cor. 4.5]{EK} that $(\zeta=0)= \mathrm{Sing} (\F)\cap {K}\neq \emptyset $ and  this implies that $\deg \left((\Omega_{{K}}^{k})^{*}\otimes\mathcal{O}_{{K}}(\deg(\F)-k)\right)>0$.
Then,
$$
\deg((\Omega_{{K}}^{k}))< \deg(\mathcal{O}_{{K}}( \deg(\F)-k)).
$$
Using (\ref{canonico}) we conclude that $\deg(\G)-k-1<  \deg(\F)-k$, i.e, $\deg(\G)\leq \deg(\F)$.

$\hfill \square$

\begin{ex}\label{hamilton} A complete flag of foliations.
\end{ex}This is an example of a complete flag of foliations to which Theorem \ref{fol1}  applies.
Let $f: \mathbb{C}^{2n} \longrightarrow \mathbb{C}$ be a polynomial function of degree $k$, write
$f = f_k + f_{k-1} + \dots + f_1$, its decomposition into homogeneous polynomials, and assume that $f$ has
only one critical point at $0 \in \mathbb{C}$. Further, suppose $(f_k =0) \subset \mathbb{P}^{2n-1}$ is a smooth algebraic
variety. The derivative of $f$ is represented by
$$
f'(z) = (\partial_1 f(z), \partial_2 f(z), \partial_3 f(z), \partial_4 f(z),\dots, \partial_{2n-1} f(z), \partial_{2n} f(z)),
$$
where $\partial_i f(z) = \frac{\partial f}{\partial z_i}(z)$. From $f'$ we can produce $2n-1$ hamiltonian vector fields
$H_i$ given by
$$\begin{array}{c}
H_1 = (-\partial_2 f(z), \partial_1 f(z), -\partial_4 f(z), \partial_3 f(z), \dots, -\partial_{2n} f(z), \partial_{2n-1} f(z)) \\
 \\
H_2= (-\partial_3 f(z), \partial_4 f(z), -\partial_1 f(z), \partial_2 f(z), \dots, -\partial_{2n} f(z), \partial_{2n-1} f(z))\\
\\
H_3= (-\partial_4 f(z), \partial_3 f(z), -\partial_2 f(z), \partial_1 f(z), \dots, -\partial_{2n} f(z), \partial_{2n-1} f(z))\\
\\
\vdots\\
\\
H_{2n-1} = (-\partial_{2n} f(z), \partial_{2n-1} f(z), \dots, -\partial_{2} f(z), \partial_{1} f(z))
\end{array}
$$
which correspond to solutions of
$$
X_1 \partial_1 f(z) + X_2 \partial_2 f(z) + \dots + X_{2n-1} \partial_{2n-1} f(z) + X_{2n} \partial_{2n} f(z)=0.
$$
Each $H_i$ is tangent to the levels $f=c$, $c \in \mathbb{C}$, and they satisfy $[H_i, H_j] =0$, $1 \leq i,j \leq 2n-1$.
Let $F(z_0, z)= f_k(z) + z_0 f_{k-1} (z) + \dots + z_0^{k-1} f_1(z)$ be the homogenized of $f$ and consider the reduced codimension one foliation $\mathcal{G}$ on $\mathbb{P}^{2n}$, of degree $k-1$, defined by
\begin{equation}\label{kupka max}
\omega = z_0 dF - k F dz_0.
\end{equation}
$\mathcal{G}$ leaves invariant the levels $f=c$ and the hyperplane at infinity $(z_0=0)= \mathbb{P}^{2n-1} \subset \mathbb{P}^{2n}$. Moreover, $\mathrm{Sing}(\mathcal{G})$ has a Baum-Kupka component $K = \{f_k=0\} \subset \mathbb{P}^{2n-1}$.

On the other hand, the vector fields $H_i$ induce a flag $\mathscr{F}$ of reduced foliations of degree $k-1$, all leaving invariant the levels $f=c$ and the hyperplane at infinity $(z_0=0)$ as follows: $\mathcal{F}_j$ is defined by $\{H_1, \dots, H_j\}$, $\dim \mathcal{F}_j =j$ and $\mathscr{F}=(\mathcal{F}_1, \mathcal{F}_2, \dots, \mathcal{F}_{2n-2}, \mathcal{F}_{2n-1} = \mathcal{G})$. Also, $\mathrm{Sing} (\F_{j+1})$ has a Baum-Kupka component $K_{j+1} \subset K$ with $K_{j+1} \not\subset \mathrm{Sing} (\F_{j})$ for $j= 1, 2, \dots, 2n-2$.

\bigskip

\noindent{\footnotesize
\textsc{Acknowlegments.} We are grateful to the referee for pointing out corrections in previous versions of this paper. The second named author is grateful to the University of Valladolid for hospitality.}

\begin{center}

\bigskip

\begin{tabular}{lll}
  Maur\'icio Corr\^ea Jr. & \hskip 20pt & M\'arcio G. Soares \\
  Dep. Matem\'atica - UFV & \hskip 20pt & Dep.Matem\'atica - UFMG \\
  Av. PH. Rolfs sn/ & \hskip 20pt & Av. Antonio Carlos 6627 \\
  36 571-000 Viçosa, Brasil & \hskip 20pt & 31 270-901 Belo Horizonte, Brasil \\
  mauricio.correa@ufv.br & \hskip 20pt & msoares@mat.ufmg.br  \\
\end{tabular}
\end{center}

\end{document}